\documentclass[11pt]{article}
\usepackage{amsmath,amssymb,a4wide}
\usepackage[curve,matrix,arrow,cmtip]{xy}
 \def\Aut{\mathop{\rm Aut}\nolimits}

 \def\Gal{\mathop{\rm Gal}\nolimits}
 \def\End{\mathop{\rm End}\nolimits}

 \def\Spec{\mathop{\rm Spec}\nolimits}
 
 \def\deg{\mathop{\rm deg}\nolimits}

\def\diag{\mathop{\rm diag}\nolimits}

\def\GL{\mathop{\rm GL}\nolimits}
\def\SL{\mathop{\rm SL}\nolimits}

\def\PGL{\mathop{\rm PGL}\nolimits}
\def\PSL{\mathop{\rm PSL}\nolimits}

\def\der{{\rm der}}

\def\sep{{\rm sep}}

\let\phi\varphi
\let\epsilon\varepsilon

\newtheorem{Thm}{Theorem}
\newtheorem{Prop}[Thm]{Proposition}

\def\qed{{\hskip0pt\unskip\unskip\nobreak\hfil\penalty50
          \hskip1em\hbox{}\nobreak\hfil
           {$\square$}
          \parfillskip=0pt\finalhyphendemerits=0
          \par}\medskip}

\newenvironment{Proof}
               {\noindent{\bf Proof.}\ }
               {\qed}

\newenvironment{Proofof}[1]
               {\noindent{\bf Proof of #1.}\ }
               {\qed}

\newcommand{\BC}{{\mathbb{C}}}

\newcommand{\BF}{{\mathbb{F}}}
\newcommand{\BG}{{\mathbb{G}}}

\newcommand{\BN}{{\mathbb{N}}}

\newcommand{\BP}{{\mathbb{P}}}
\newcommand{\BQ}{{\mathbb{Q}}}

\newcommand{\BZ}{{\mathbb{Z}}}

\newcommand{\CF}{{\cal F}}

\newbox\mybox
\def\arrover#1{\mathrel{
       \setbox\mybox=\hbox spread 1.4em
              {\hfil$\scriptstyle#1$\hfil}
       \vbox{\offinterlineskip\copy\mybox
             \hbox to\wd\mybox{\rightarrowfill}}}}

\def\larrover#1{\mathrel{
       \setbox\mybox=\hbox spread 1.4em
              {\hfil$\scriptstyle#1\vphantom{g}$\hfil}
       \vbox{\offinterlineskip\copy\mybox
             \hbox to\wd\mybox{\leftarrowfill}}}}

\def\ontoover#1{\mathrel{
       \setbox\mybox=\hbox spread 1.4em
              {\hfil$\scriptstyle#1\vphantom{g}$\hfil}
       \vbox{\offinterlineskip\copy\mybox
             \hbox to\wd\mybox{\rightarrowfill\hskip-2.8mm
                               $\rightarrow$}}}}
\def\leftontoover#1{\mathrel{
       \setbox\mybox=\hbox spread 1.4em
              {\hfil$\scriptstyle#1\vphantom{g}$\hfil}
       \vbox{\offinterlineskip\copy\mybox
             \hbox to\wd\mybox{$\leftarrow$\hskip-2.8mm
                               \leftarrowfill}}}}
\let\longto\longrightarrow
\let\into\hookrightarrow

\def\Cinf{{\BC}_\infty}

\def\kinf{k_\infty}

\renewcommand{\Im}{{\mathrm {Im}}}

\begin{document}

\title{Galois groups associated to generic Drinfeld modules and a conjecture of Abhyankar \\ {\bf Notice of Replacement}}
\author{Florian Breuer\\ 
Stellenbosch University, Stellenbosch, South Africa\\ 
{\em fbreuer@sun.ac.za}}
\maketitle


The previous version of this paper, as disseminated on arXiv on 17 March 2013, contains an important gap at the end of the proof. It is superceded by the paper
\cite{new}. 

The erroneous claim is that $K_T = K_{T,0}(\lambda(v))$ is a transcendental extension of $K_{T,0}$, whereas in fact $\lambda(v)$ satisfies the equation
\[
\lambda(v)^{q^r-1} = T\prod_{0\neq w\in V}\frac{\lambda(v)}{\lambda(w)} \in K_{T,0}.
\]
To correct this error, one effectively needs to show that $K_T\cap K_{NT,0} = K_{T,0}$ inside $K_{NT}$. This is achieved in \cite{new} by studying the fixed field of $G(NT,T)\cap\SL_r(A/NTA)$ inside $K_{TN,0}$ via Anderson's determinant morphism from $M_{TN}^r$ to $M_{TN}^1$, and applying some group theory.

The paper \cite{new} also contains an explicit construction of $M_{TN}^r$, extending the construction of Richard Pink, which may be of independent interest.

Many thanks to an alert anonymous referee at the {\em Journal of Number Theory} for pointing out this error.

\end{document}